\documentclass [winedt,yap]{iitparc}
\usepackage{cite}
\usepackage{amsmath,amssymb,amsfonts,bm}
\usepackage{lscape}
\usepackage{floatfig,wrapfig,epsfig}
\usepackage{subfigure}
\usepackage{color}
\usepackage{psboxit}
\usepackage{rotating}
\usepackage{curves}
\usepackage{ulem}
\usepackage[mathscr]{eucal}
\RequirePackage{psfrag}
\RequirePackage{graphicx}

\begin{document}\normalem
%
\frontmatter          
%
%
\IssuePrice{25.00}%
\TransYearOfIssue{2006}%
\TransCopyrightYear{2006}%
\OrigYearOfIssue{2006}%
\OrigCopyrightYear{2006}%
\TransVolumeNo{67}%
\TransIssueNo{3}%
\OrigIssueNo{3}%
%
\mainmatter              
%
\setcounter{page}{480}
\CRubrika{CONTROL IN SOCIAL ECONOMIC SYSTEMS}
\Rubrika{CONTROL IN SOCIAL ECONOMIC SYSTEMS}
\newtheorem*{conj*}{Conjecture}{\bfseries}{\itshape}

\title{Analytical Expression of the Expected Values of Capital \\ at Voting
in the Stochastic Environment}

\titlerunning{Expected Values of Capital at Voting}
\author{P.  Yu.  Chebotarev}

\authorrunning{Chebotarev}

\OrigCopyrightedAuthors{Chebotarev}

\institute{Trapeznikov Institute of Control Sciences, Russian Academy of
Sciences, Moscow, Russia}

\received{Received February 01,  2005}

\OrigPages{152--165}

\maketitle

\begin{abstract}
In the simplest variant of the model of collective decision making in the
stochastic environment, the participants were segregated into egoists and a
group of collectivists.  ``Proposal of the environment'' is the
stochastically generated vector of algebraic increments of capitals.  The
social dynamics was defined by the sequence of proposals accepted by
threshold-majority voting.  Analytical expressions of the expected values
of the capitals of participants, collectivists and egoists were obtained.
Distinctions of some principles of group voting were discussed.
\end{abstract}

\DOI{012X}

\section{introduction}

A model where the participants vote for the projects of redistribution of
their own property was analyzed by A.V.~Malishevskii in the late 1960's [1,
pp.~93--95].  In this model, voting is greatly manipulatable by the
organizers, as it is the case with the participants whose ``ideals'' are
the points in the multidimensional space of programs [2] (see also [3]).
The monograph [4] is devoted to the spatial voting models.  The problems of
relations between selfishness, altruism, and rationality were considered in
[5--7], and voting as a method of making decisions about redistribution of
social benefits by means of taxation and social programs was discussed in
[8--11].

For a stochastic environment oriented to the analysis of ``effectiveness''
of the voters' collectivist and selfish attitudes in the conditions where
new programs are randomly generated by the ``environment,'' rather than
developed by the organizers or participants of voting, a model of voting
was suggested by the present author and analyzed in [12].  Consideration
was given to the case of neutral environment, the rest of the cases being
analyzed mostly in qualitative terms.  Emphasis was made on the time
dependencies of the participant ``capitals'' under various values of the
model parameters.  In what follows, the present author obtained explicit
expressions of these dependencies, the exact formulas of expectations and
their normal approximations.  These formulas enable one to determine the
nature of social dynamics under any values of the model parameters.
Expressions for the expectations of the capital increments were obtained in
Sec.~\ref{sec_main}, the necessary lemmas being proved in
Sec.~\ref{sec_lem}.

We give a thumbnail of the simplest variant of the model.  The ``society''
is assumed to consist of $n$ ``participants'' of which $\ell$ are
``egoists'' and $g=n-\ell$ are the ``group'' members.  At any time instant,
the participant is characterized by its ``capital'' expressed by a real
number and interpreted in the most general sense (like {\em utility}).
Each participant has some starting capital.  The ``environment proposal''
is the vector $(d_1,\ldots,d_n)$ of the algebraic increments of the
capitals of all participants.  These increments are
independent-in-aggregate random variables with identical distribution
$N(\mu,\sigma^2)$, where $\mu$ and $\sigma$ are the model parameters.  At
each step, one random ``environment proposal'' is put to the vote.  Each
``egoist'' votes ``for'' if and only if the proposal brings it a positive
capital increment.  Each member of the group votes ``for'' if and only if
the group gains from the realization of this proposal.  The ``gain'' may be
understood differently.  The basic model considers two fundamental
principles of decision voting .

{\bf Principle A}.  {\it The group votes ``for'' the proposal $(d_1,\ldots,
d_n)$ if and only if as the result of its approval the number of group
members getting a positive capital increment exceeds that of the group
members getting a negative increment.}

{\bf Principle B}.  {\it The group votes ``for'' the proposal $(d_1,\ldots,
d_n)$ if and only if the sum of increments of group members is positive$:$
$\sum d_i>0$, where the sum is taken over the subscripts of the group
participants.}

The results of voting are summarized using the ``$\alpha$-majority''
procedure:  the proposal is accepted if and only if more than $\alpha n$,\
$0\leq\alpha<1$, participants vote for it.  If the proposal is
accepted, then the capitals of all participants get the corresponding
increments $d_1,\ldots,d_n$, or remain the same, otherwise.

Capital dynamics of the participants is analyzed in terms of their social
roles (``egoist'' or ``group member'') and model parameters.  It is
implied, in particular, that a scenario is plausible where the egoists join
the group and the ``group egoism'' resembles more and more the decision
making in the interests of the entire society.  By another hypothetical
scenario, an ``ineffective'' group dissolves, and its members either become
``egoists'' or make other groups.  It is planned also to consider the case
of a socially oriented group supporting its poorest members and preventing
their ruin, a variant of the model where the capital increments depend on
their current values, the impact on the social dynamics of the mechanisms
of taxation and collection of the ``party dues,'' and so on.  At the same
time, it is not planned to consider purely economical mechanisms of capital
reproduction and loss because our aim lies in analyzing the social, rather
than the economic phenomena.

Consideration is given not only to the traditional decision threshold
$\alpha=0{.}5$ corresponding to the ``simple majority,'' but to all
thresholds ranging from $0$ to $1$.  That is due to the fact
that the most important---for example, ``constitutional''---decisions are
accepted by the a ``qualified parliamentary majority'' with a threshold
greater than $0{.}5$.  On the other hand, there are ``initiative''
decisions such as forming new deputy groups in parliament, putting question
on agenda, sending requests to other state authorities, initiating
referendums, and so on which can be approved by a certain number of votes
smaller than one half.  This model is discussed in more detail in [12].

The relation between the numbers of egoists and group members is defined by
the parameter $\beta=\ell/2n$, half of the portion of egoists among all
participants.  We denote by $f_{\mu,\sigma}(\cdot)$ and $F_{\mu,\sigma}(
\cdot)$, respectively, the one-dimensional density and the cumulative
distribution function corresponding to the distribution $N(\mu,\sigma^2)$;\
$f(\cdot)$ and $F(\cdot)$ stand for density and the distribution function
of the normal distribution with center at $0$ and variance $1$;\ $M(\xi)$
and $\sigma(\xi)$ are, respectively, the expectation and the deviation
of any of the random variable $\xi$ at hand.

Each $i$th egoist participant votes for a proposal if and only if $d_i>0$.
The probability of this event is as follows:
\begin{gather}\label{p}
p=P\{d_i>0\}=1-F_{\mu,\sigma}(0)=F\left(\frac{\mu}{\sigma}\right);
\end{gather}
the probability of voting ``against'' the proposal is as follows:
\begin{gather}\label{q}
q=1-p=P\{d_i<0\}=F_{\mu,\sigma}(0)=1-F\left(\frac{\mu}{\sigma}\right)=F\left(-\frac{\mu}{\sigma}\right).
\end{gather}
In what follows, we also need the notation
\begin{gather}\label{f}
f=f\left(\frac{\mu}{\sigma}\right).
\end{gather}

According to the model, the probability that a participant abstains from
voting is zero because the normal distribution is continuous.  Therefore,
the voting of an egoist participant is the Bernoulli test with the
parameter $p$.  Then, since the values $d_i$ are independent, the number of
egoists voting ``for'' is distributed {\it binomially\/} with the
parameters $\ell$ and $p$.  The mean value and the variance of this
distribution are, respectively, $p\ell$ and $pq\ell$.

\section{lemmas of the ``normal voting sample''} \label{sec_lem}

In this section we prove lemmas that underlie the following calculations.
By the ``{\it normal voting sample}'' of size%
\footnote[1]{In this section, $\ell$ is an arbitrary integer.}
$\ell$ with the parameters $(\mu,\sigma^2)$ and voting threshold $\ell_0$
is meant the totality of random variables $\left(\zeta_1I(\overline{\zeta},
\ell_0),\ldots,\zeta_{\ell}I(\overline{\zeta},\ell_0)\right)$, where
${\overline{\zeta}=(\zeta_{1},\ldots,\zeta_{\ell})}$ is a sample from the
distribution $N(\mu,\sigma^2)$,
\begin{gather}\label{vartheta}
I(\overline{\zeta},\ell_0)
=\left\{%
\begin{array}{ll}
    1& \text{if\ } n^+(\overline{\zeta})>\ell_0;\\
    0, & \text{otherwise,}
\end{array}
\right.
\end{gather}
and%
\footnote[2]{$\#X$ denotes the number of elements in the finite set $X$.}
\begin{gather}\label{n^+}
n^+(\overline{\zeta})=\#\{k:\zeta_k>0,\,k=1,\ldots,\ell\}.
\end{gather}

The following lemma holds.

\begin{lemma}\label{l1}
Let $(\eta_1,\ldots,\eta_{\ell})$ be a normal voting sample with the
parameters $(\mu,\sigma^2)$ and the voting threshold $\ell_0$.  Then,%
\footnote[3]{Here and below the sum is zero if the lower limit is greater
than the upper one.  The integer part is bracketed.}
for any $k=1,\ldots,\ell$
\begin{gather}\label{Mvot}
M(\eta_k) =\sum_{x=[\ell_0]+1}^{\ell}\left(\mu+\frac{\sigma
f}{q}\left(\frac{x}{p\ell}-1\right)\right)
 \left(\begin{array}{c} \ell \\ x  \end{array}\right)
p^xq^{\ell-x},
\end{gather}
where $p,\,q$, and $f$ are defined in \ntt{\ref{p}}{\ref{f}}.
\end{lemma}

The proofs are given in the Appendix.

To calculate (\ref{Mvot}), one may use the well-known relation between the
binomial distribution and the beta-distribution:
\begin{gather}\label{BetaD}
\sum_{x=t}^{\ell} \left(\begin{array}{c} \ell \\ x
\end{array}\right)p^xq^{\ell-x} =B(p\mid t,\ell-t+1),
\end{gather}
the right-hand side contains the cumulative function of the distribution of
beta-distribution with $t$ and $\ell-t+1$ degrees of freedom.  However, the
normal approximation of the binomial probability can be used even for a
comparatively small $\ell$:
\begin{gather}\label{appro0}
 \left(\begin{array}{c} \ell \\ x  \end{array}\right)p^xq^{\ell-x}
 \approx f\left(\frac{x-\mu'}{\sigma'}\right) \approx
                                 F\left(\frac{x+0{.}5-\mu'}{\sigma'}\right) -
                                 F\left(\frac{x-0{.}5-\mu'}{\sigma'}\right),
\end{gather}
where $\mu'$ and $\sigma'$ are the same parameters as for the binomial
distribution:  $\mu'=p\ell$, $\sigma' =\sqrt{pq\ell}$.  Summation in
(\ref{appro0}) provides an approximation to the binomial distribution
function:
\begin{gather}\label{appro}
\begin{split}
 &\sum_{x=0  }^{ t} \left(\begin{array}{c} \ell \\ x
\end{array}\right)p^xq^{\ell-x}\approx F\left(
\frac{t+0{.}5-p\ell}{\sqrt{pq\ell}}\right) \quad \text{and}
\\
&\sum_{x=t+1}^{\ell} \left(\begin{array}{c} \ell \\ x
\end{array}\right)p^xq^{\ell-x}\approx
F\left(-\frac{t+0{.}5-p\ell}{\sqrt{pq\ell}}\right).
\end{split}
\end{gather}

The normal approximation is advisable for $pq\ell\geq9$.  For a fixed
$pq\ell$, its accuracy is maximal for $p=0{.}5$ and decreases with $p$
approaching $0$ or $1$.  Therefore, the normal approximation is often
recommendable for $0{.}1<p<0{.}9$ already if $pq\ell>5$.  For the values of
$p$ that are very close to $0$ or $1$, usually the condition $pq\ell>25$ is
imposed.

\begin{lemma}\label{l2}
Let $(\eta_1,\ldots,\eta_{\ell})$ be a normal voting sample with the
parameters $(\mu,\sigma^2)$ and the voting threshold $\ell_0$.  Then
substitution of the standard normal approximation of the random variable
$n^+(\overline{\zeta})$ for its binomial distribution provides the
following approximation of $M(\eta_k)$ for any $k=1,\ldots,\ell:$
\begin{gather}\label{MvotI}
M(\eta_k) \approx\mu F(-\ell'_0)+\frac{\sigma f}{\sqrt{
pq\ell}}f(\ell'_0),
\end{gather}
where $\ell'_0=([\ell_0]+0{.}5-p\ell)/\sqrt{pq\ell}$.
\end{lemma}

This approximation is equivalent to the replacement of a finite number of
the participants of voting by a continuous ``voting field.'' We present for
completeness sake the following lemma where $J(\overline{\zeta}, \ell_0)$
differs from $I(\overline{\zeta},\ell_0)$ (\ref{vartheta}) in the sign of
inequality in the definition.  Its proof is similar to that of
Lemmas~\ref{l1} and~\ref{l2}.

\begin{lemma}\label{l3}
Let $(\eta_1,\ldots,\eta_{\ell})=(\zeta_1 J(\overline{\zeta},\ell_0),
\ldots,\zeta_{\ell} J(\overline{\zeta},\ell_0))$, where $\overline{\zeta}=
(\zeta_{1},\ldots,\zeta_{\ell})$, be a sample from the distribution
$N(\mu,\sigma^2)$ and
\begin{gather}\label{vartheta1}
J(\overline{\zeta},\ell_0)
=\left\{%
\begin{array}{ll}
    1& \text{if}\quad n^+(\overline{\zeta})\leq\ell_0;\\
    0, & \text{otherwise.}
\end{array}
\right.
\end{gather}
Then, for any $k=1,\ldots,\ell$
\begin{gather}\label{Mivot}
M(\eta_k) =\sum_{x=0}^{[\ell_0]}\left(\mu+\frac{\sigma
f}{q}\left(\frac{x}{p\ell}-1\right)\right) \left(\begin{array}{c}
\ell \\ x
\end{array}\right)p^xq^{\ell-x},
\end{gather}
and the normal approximation provides the approximation
\begin{gather}\label{MvotI'}
M(\eta_k)\approx\mu F(\ell'_0)-\frac{\sigma f}{\sqrt{
pq\ell}}f(\ell'_0),
\end{gather}
where $\ell'_0=([\ell_0]+0{.}5-p\ell)/\sqrt{pq\ell}$.
\end{lemma}

To sum up the auxiliary results, we denote by $\mu^+(\mu,\sigma,\ell,
\ell_0)$ the expectation of the normal voting sample of size $\ell$ with
the parameters $(\mu,\sigma^2)$ and the voting threshold $\ell_0$, and by
$\mu^-(\mu,\sigma,\ell,\ell_0)$, the expectation of the set of random
variables as defined in Lemma~\ref{l3}.  Then, the exact and approximate
values of $\mu^+(\mu,\sigma,\ell,\ell_0)$ obey the formulas (\ref{Mvot})
and (\ref{MvotI}), and the exact and approximate values of
$\mu^-(\mu,\sigma,\ell,\ell_0)$ obey the formulas (\ref{Mivot}) and
(\ref{MvotI'}).

\section{general expressions for the mean increments of capitals}
\label{sec_main}

\subsection{Increment of the Egoist's Capital}

Let us calculate the expectation $M(\widetilde{d}_{\mathcal{E}})$ of the
egoist's capital increment in one step.%
\footnote[4]{The notation under tilde refers to the actual, in contrast to
the proposed, capital increments.}
We recall that a step lies is considering one proposal of the environment
by the participants, independently of whether it is accepted or not.  Let
the event $G$ be supported by the group, and $\overline{G}$ be the
opposite. Then,
\begin{gather}\label{MDe}
M(\widetilde{d}_{\mathcal{E}})= M(\widetilde{d}_{\mathcal{E}}\mid
G)P\{G\}+M(\widetilde{d}_{\mathcal{E}}\mid\overline{G})P\{\overline{G}\}.
\end{gather}
We denote $P\{G\}$ and $P\{\overline{G}\}=1-P\{G\}$, respectively, by $P_G$
and $Q_G$.

Support of a proposal by the group brings it $(1-2\beta)n$ votes.
Consequently, for approval of a proposal on condition that it is supported
by the group, it is necessary and sufficient that more than
$(\alpha-(1-2\beta))n$ egoists vote for it.  If the group does not support
the proposal, then it is necessary and sufficient that more than $\alpha n$
egoists vote for it.  Therefore,
\begin{gather}\label{MDeG}
M(\widetilde{d}_{\mathcal{E}}\mid G)=\mu^+(\mu,\sigma,\ell,\gamma
n),\quad
M(\widetilde{d}_{\mathcal{E}}\mid\overline{G})=\mu^+(\mu,\sigma,\ell,\alpha
n),
\end{gather}
where $\gamma=\alpha-(1-2\beta)$ and $[t]$ is the integer part of $t$.
Hence,
\begin{gather}\label{MDe'}
M(\widetilde{d}_{\mathcal{E}}) =\mu^+(\mu,\sigma,\ell,\gamma n)P_G
+\mu^+(\mu,\sigma,\ell,\alpha n)Q_G.
\end{gather}

By substituting expression (\ref{Mvot}) for $\mu^+(\cdot,\cdot,\cdot,\cdot)
$ in (\ref{MDe'}) and using the notation
\begin{gather}\label{bxl}
b(x\mid\ell)= \left(\begin{array}{c} \ell \\ x
\end{array}\right)p^xq^{\ell-x},
\end{gather}
we obtain
\begin{eqnarray}\label{MDe''}
M(\widetilde{d}_{\mathcal{E}})
&=&P_G\sum_{x=[\gamma n]+1}^{\ell}\left(\mu+\frac{\sigma
 f}{q}\left(\frac{x}{p\ell}-1\right)\right)b(x\mid\ell)
 \\
&+&Q_G   \sum_{x=[\alpha    n]+1}^{ \ell} \left(\mu+\frac{\sigma
f}{q}\left(\frac{x}{p\ell}-1\right) \right)b(x\mid\ell)
\nonumber
\\
\nonumber
&=& P_G \sum_{x=[\gamma n]+1}^{[\alpha n]} \left(\mu+\frac{\sigma f}{q}
\left(\frac{x}{p\ell}-1\right) \right)b(x\mid\ell)
\nonumber
\\
&+&           \sum_{x=[\alpha    n]+1}^{ \ell}
\left(\mu+\frac{\sigma f}{q}\left(\frac{x}{p\ell}-1\right)
\right)b(x\mid\ell).\nonumber
\end{eqnarray}
If the group follows Principle A, then
\begin{gather}\label{PGA}
P_G^A=\sum_{x=[{g}/{2}]+1}^gb(x\mid g)\quad\text{and}\quad
Q_G^A=1-P_G^A
     =\sum_{x=0}^{[{g}/{2}]}b(x\mid g)
\end{gather}
should be substituted for $P_G$ and $Q_G$, respectively, and in the case of
Principle~B,
\begin{gather}\label{PGB}
P_G^B=F\left(\frac{\mu\sqrt{g}}{\sigma}\right)\quad\text{and}\quad
Q_G^B=1-P_G^B
     =F\left(-\frac{\mu\sqrt{g}}{\sigma}\right),
\end{gather}
because for an arbitrary proposal the mean capital increment of the group
has the distribution $N\left(\mu,\dfrac{\sigma}{\sqrt{g}}\right)$.

Now, by substituting in (\ref{MDe'}) the normal approximation (\ref{MvotI})
for $\mu^+(\cdot,\cdot,\cdot,\cdot)$ and using the notation
\begin{gather}\label{Fth}
F_{\theta}=F\left(-\frac{[\theta
n]+0{.}5-p\ell}{\sqrt{pq\ell}}\right),\quad f_{\theta}=f\left(
\frac{[\theta n]+0{.}5-p\ell}{\sqrt{pq\ell}}\right),
\end{gather}
we get
\begin{eqnarray}\label{MDeAp}
M(\widetilde{d}_{\mathcal{E}}) 
&\approx&
 P_G\left(\mu F_{\gamma}+\frac{\sigma f}{\sqrt{ pq\ell}}f_{\gamma}\right)
+Q_G\left(\mu F_{\alpha}+\frac{\sigma f}{\sqrt{
pq\ell}}f_{\alpha}\right)
\\
\nonumber 
&=&[P_G\;Q_G]\left(\mu
\left[\begin{array}{c} F_{\gamma}\\
F_{\alpha}
\end{array}
\right]
                +\frac{\sigma f}{\sqrt{ pq\ell}}\left[\begin{array}{c} f_{\gamma}\\
f_{\alpha}
\end{array}
\right]\right)
\end{eqnarray}
where matrix notation is used.  In the case of Principle~A, the approximate
values
\begin{gather}\label{PGAQGA}
\begin{split}
  &F\left(-\frac{\left[\tfrac{g}{2}\right]+0{.}5-pg}{\sqrt{pqg}}\right)
\approx P_G^A,
 \\
&1-F\left(-\frac{\left[\tfrac{g}{2}\right]+0{.}5-pg}{\sqrt{pqg}}\right)=
  F\left( \frac{\left[\tfrac{g}{2}\right]+0{.}5-pg}{\sqrt{pqg}}\right)
\approx Q_G^A
\end{split}
\end{gather}
or exact values (\ref{PGA}) must be substituted in this formula; in the
case of Principle~B, the exact values of (\ref{PGB}).

\subsection{Increment of the Capital of Group Member}

Now we derive formulas similar to (\ref{MDe''}) and (\ref{MDeAp}) for the
expectation $M(\widetilde{d}_{\mathcal{G}})$ of the capital increment of a
group member in one step.  Let $E_{\theta}$ be an event of egoists giving
more than $\theta n$ votes to the proposal; $\overline{E}_{\theta}$,
otherwise.  If $E_{\alpha}$ is satisfied, then the proposal is accepted
independently of the group voting; if $\overline{E}_{\alpha}\wedge
E_{\gamma}$ takes place, then acceptance/rejection of the proposal is
defined by its support by the group.  If $\overline{E}_{\gamma}$ is
realized, then the proposal cannot be accepted.  Therefore,
\begin{gather}\label{MdG}
M(\widetilde{d}_{\mathcal{G}}) =M(\widetilde{d}_{\mathcal{G}}\mid
E_{\alpha}) \,    P\{E_{\alpha}\}+
        M(\widetilde{d}_{\mathcal{G}}\mid\overline{E}_{\alpha}\wedge E_{\gamma})\,
        P\{\overline{E}_{\alpha}\wedge E_{\gamma}\}
       +0                                          \cdot P\{\overline{E}_{\gamma}\}.
\end{gather}
If $E_{\alpha}$ is satisfied, then the proposal is always accepted.
Therefore, ${M(\widetilde{d}_{\mathcal{G}}\!\mid\!  E_{\alpha})\!=\!\mu}$
by virtue of independence of the capital increments.  If
${\overline{E}_{\alpha}\!\wedge\!  E_{\gamma}}$ is satisfied, then
acceptance of the proposal is defined by the group voting, and generally
${M(\widetilde{d}_{\mathcal{G}}\!\mid\!\overline{E}_{\alpha}\!\wedge\!
E_{\gamma})\!\ne\!\mu}$.

We make use of the notation
\begin{gather}\label{Pthe}
P_{\theta} =P\{E_{\theta}\} =\sum_{x=[\theta
n]+1}^{\ell}b(x\mid\ell) \approx F_{\theta}=F\left(-\frac{[\theta
n]+0{.}5-p\ell}{\sqrt{pq\ell}}\right),
\end{gather}
to rearrange (\ref{MdG}) in
\begin{gather}\label{MdG'}
M(\widetilde{d}_{\mathcal{G}}) =\mu P_{\alpha}
+M(\widetilde{d}_{\mathcal{G}}\mid\overline{E}_{\alpha}\wedge
E_{\gamma})(P_{\gamma}-P_{\alpha}).
\end{gather}

If the group follows Principle~A, then $M(\widetilde{d}_{\mathcal{G}}\mid
\overline{E}_{\alpha}\wedge E_{\gamma})=\mu^+(\mu,\sigma,g,g/2)$, the
substitution of (\ref{Mvot}) in
\begin{gather}\label{MdGA}
M(\widetilde{d}_{\mathcal{G}}^A) =\mu P_{\alpha}
+\mu^+\left(\mu,\sigma,g,g/2\right)(P_{\gamma}-P_{\alpha}),
\end{gather}
provides
\begin{eqnarray}\label{MdG''}
M(\widetilde{d}_{\mathcal{G}}^A) 
&=& \mu P_{\alpha}
+\sum_{x=[g/2]+1}^{g}\left(\mu+\frac{\sigma
f}{q}\left(\frac{x}{pg}-1\right) \right)b(x\mid
g)(P_{\gamma}-P_{\alpha})
\\
\nonumber
&=&
 \mu \sum_{x=[\alpha n]+1}^{\ell}b(x\mid\ell)
+       \sum_{x=[g/2]+1}^{g}\left(\mu+\frac{\sigma
f}{q}\left(\frac{x}{pg}-1\right) \right)b(x\mid g)
     \sum_{x=[\gamma n]+1}^{[\alpha n]}b(x\mid\ell).
\end{eqnarray}

In order to obtain the normal approximation of this value, we replace
$\mu^+(\mu,\sigma,g,g/2)$ in (\ref{MdGA}) by the approximation
(\ref{MvotI}), and $P_{\alpha}$ and $P_{\gamma}$, by the approximations
$F_{\alpha}$ and $F_{\gamma}$ (see (\ref{Pthe})).  By using the notation
\begin{gather}\label{fGA}
f_G^A
=f\left(\frac{\left[\tfrac{g}{2}\right]+0{.}5-pg}{\sqrt{pqg}}\right),
\end{gather}
we obtain
\begin{eqnarray}\label{MdGAp}
M(\widetilde{d}_{\mathcal{G}}^A)  
&\approx& \mu F_{\alpha}+\left(\mu
P_G^A+\frac{\sigma
ff_G^A}{\sqrt{pqg}}\right)(F_{\gamma}-F_{\alpha})
\\
&=& \mu(P_G^A F_{\gamma}+(1-P_G^A)F_{\alpha})+\frac{\sigma
ff_G^A}{\sqrt{pqg}}(F_{\gamma}-F_{\alpha})
\nonumber
\\
 \nonumber
&=&[F_{\gamma}\;F_{\alpha}]\left(\mu\left[
\begin{array}{c}
P_G^A
 \\[1mm]
Q_G^A
\end{array}
\right] +\frac{\sigma ff_G^A}{\sqrt{pqg}}\left[
\begin{array}{r}
1 \\ -1
\end{array}
\right]\right).
\end{eqnarray}

If the group follows Principle~B, then ${M(\widetilde{d}_{\mathcal{G}}\!
\mid\!\overline{E}_{\alpha}\!\wedge E_{\gamma})\!=\!\mu^+\!\left(\mu,\frac{
\sigma}{\sqrt{g}},1,0\right)}$, and (\ref{MdG'}) is reshaped in
\begin{gather}\label{MdGB}
M(\widetilde{d}_{\mathcal{G}}^{\,B}) =\mu P_{\alpha}
+\mu^+\left(\mu,{\sigma\over\sqrt{g}},1,0\right)(P_{\gamma}-P_{\alpha}).
\end{gather}

To determine $\mu^+\left(\mu,\dfrac{\sigma}{\sqrt{g}},1,0\right)$, one may
substitute the corresponding arguments in (\ref{Mvot}) or, which is even
simpler, multiply (A.4) by $p$:
\begin{gather}\label{mu+B}
\mu^+\left(\mu,\frac{\sigma}{\sqrt{g}},1,0\right) =\mu
F\left(\frac{\mu\sqrt{g}}{\sigma}\right)+\frac{\sigma}{\sqrt{g}}f\left(\frac{\mu\sqrt{g}}{\sigma}\right).
\end{gather}

By using the notation (\ref{PGB}) and
\begin{gather}\label{fGB}
f_G^B =f\left(\frac{\mu\sqrt{g}}{\sigma}\right),
\end{gather}
we obtain
\begin{eqnarray}\label{MdGBAp}
M(\widetilde{d}_{\mathcal{G}}^{\,B}) 
&\approx& \mu P_{\alpha}+\left(\mu P_G^B+\frac{\sigma
f_G^B}{\sqrt{g}}\right)(P_{\gamma}-P_{\alpha})
\\
&=&[P_{\gamma}\;P_{\alpha}]\left(\mu\left[\begin{array}{c}
 P_G^B \\[1mm] Q_G^B
\end{array}
\right] +\frac{\sigma f_G^B}{\sqrt{g}}\left[
\begin{array}{r}
 1 \\ -1
\end{array}\right]\right),
 \nonumber
\end{eqnarray}
where either the exact values or normal approximations (see (\ref{Pthe}))
can be substituted for $P_{\gamma}$ and $P_{\alpha}$.  The expected values
of the capital increments of the egoists and members of the group in a
series of $s$ steps obey, respectively, $sM(\widetilde{d}_{\mathcal{E}})$
and $sM(\widetilde{d}_{\mathcal{G}})$.

\subsection{Voting Principles A$'$ and A$''$}

In addition to Principles~A and~B, the group can make use of
Principle~A$'$ [12] enabling it sometimes to minimize and sometimes
to eliminate the advantage of the egoists in those zones of variation
of the threshold $\alpha$ where they have this advantage.  This is
reached at the expense of reduced group capital increment in these
zones:  relative gain of the group leads to absolute losses.  We formulate
below this voting principle.  The threshold $\alpha'$ vs.~$\alpha$
and~$\beta$ is depicted in Fig.~1.

\begin{figure}[t]
\centerline{\includegraphics[clip]{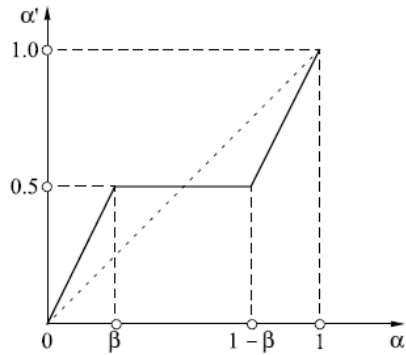}}
\caption{Threshold $\alpha'$ vs.  $\alpha$.}
\end{figure}

{\bf Principle~A$'$}.  {\it The group votes ``for'' the proposal
of the environment if and only if as the result of accepting it
the part of its members getting a positive capital increment
exceeds the threshold}
\begin{eqnarray}\label{aa'}
   \alpha'&=&\left\{%
\begin{array}{lll}
    \tfrac{1}{2}-\tfrac{\delta}{2\beta}, &\text{where\ }\delta=\beta-\alpha,    & \text{for\
    }\alpha<\beta;\\[3mm]
    \tfrac{1}{2}+\tfrac{\delta}{2\beta}, &\text{where\ }\delta=\alpha-(1-\beta),& \text{for\
    }\alpha>1-\beta;\\[3mm]
    \tfrac{1}{2}                      &                               & \text{for\
    }\alpha\in[\beta,1-\beta]
\end{array}%
\right.\hspace*{4mm} 
 \\ &=&
\left\{%
\begin{array}{ll}
    \tfrac{\alpha}{2\beta}  & \text{for\ }\alpha<\beta;\\[1mm]
    1-\tfrac{1-\alpha}{2\beta}  & \text{for\ }\alpha>1-\beta;\\[1mm]
    \tfrac{1}{2} & \text{for\ }\alpha\in[\beta,1-\beta].
\end{array}%
\right.\nonumber
\end{eqnarray}

\begin{figure}[t]
\centerline{\includegraphics[clip]{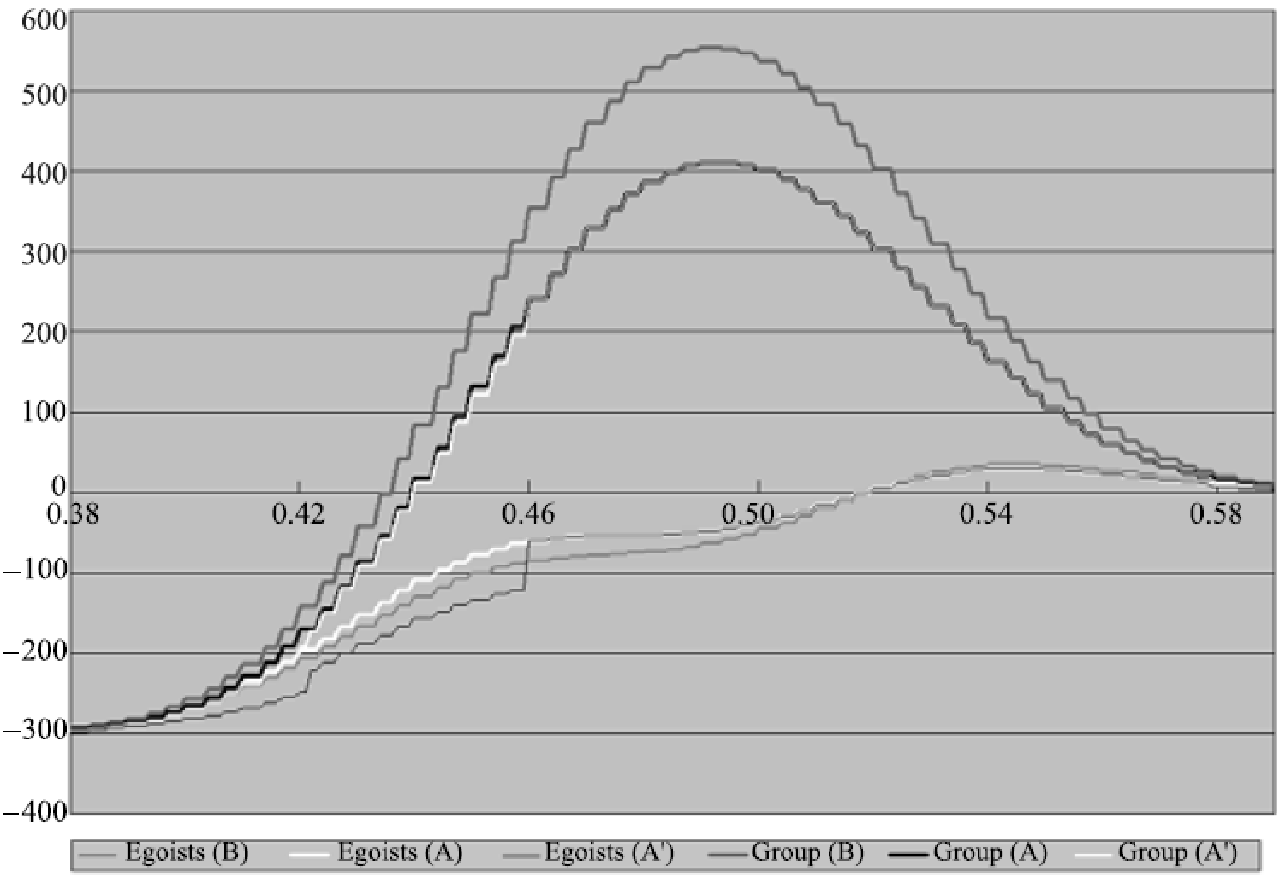}}
\caption{Capital increments of the egoists and group members following
Principles~A, B, and A$'$ vs.  the decision threshold $\alpha$.  Values of
the parameters:  $\beta=0{.}46$; $n=300$; $\mu=-0{.}3$; $\sigma=10$;
$s=1000$.}
\bigskip
\end{figure}


\begin{figure}[t]
\centerline{\includegraphics[clip]{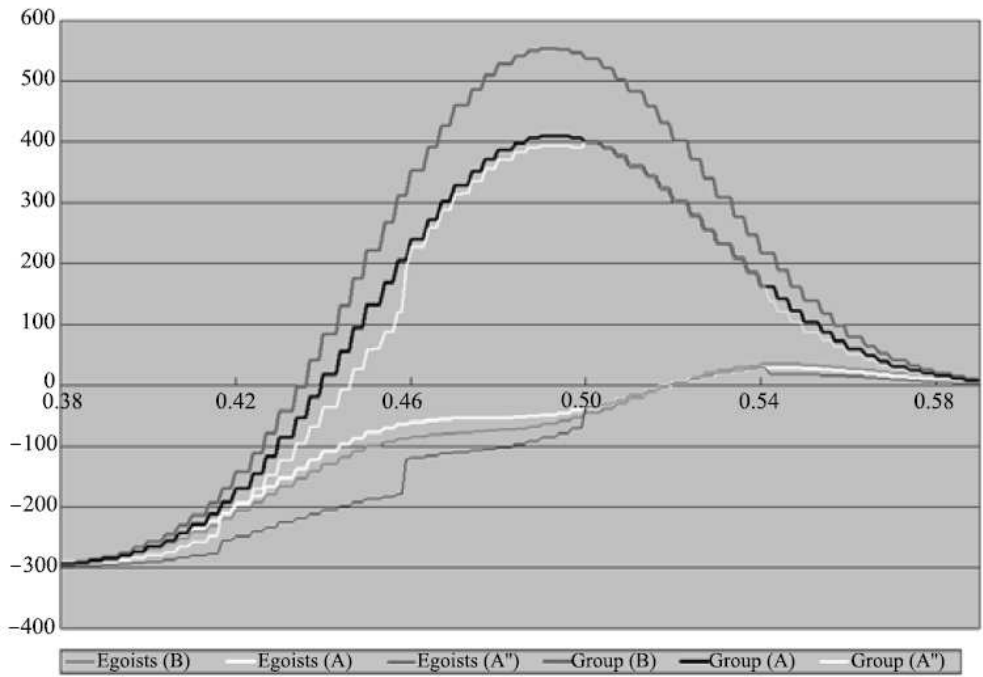}}
\caption{Capital increments of egoists and group members following
Principles~A, B, and~A$''$ vs.  the decision threshold $\alpha$.  Values of
the parameters:  $\beta=0{.}46$; $n=300$; $\mu=-0{.}3$; $\sigma=10$;
$s=1000$.}
\bigskip
\end{figure}

Since Principle~A$'$ differs from Principle~A only in the internal
threshold of group voting, for calculation of the expected capital
increments it suffices to replace $\left[{g}/{2} \right]$ by
$[\alpha'g]$ in (\ref{PGA}), (\ref{PGAQGA}), (\ref{MdG''}), and
(\ref{fGA}).  Thus, in (\ref{MDe''}), (\ref{MDeAp}), and (\ref{MdGAp}) one
has to substitute
\begin{gather}\label{PGA'}
\begin{split}
&P_G^{A'}=\sum_{x=[\alpha'g]+1}^gb(x\mid g)
\\
& \text{or the approximation}\quad
F\left(-\frac{[\alpha'g]+0{.}5-pg}{\sqrt{pqg}}\right)\approx
P_G^{A'},
\end{split}
\end{gather}
for $P_G$ and $P_G^A$,
\begin{gather}\label{QGA'}
\begin{split}
&Q_G^{A'}=1-P_G^{A'}
     =\sum_{x=0}^{[\alpha'g]}b(x\mid g)
\\
  &    \text{or the  approximation}
      \quad
F\left( \frac{[\alpha'g]+0{.}5-pg}{\sqrt{pqg}}\right)\approx
Q_G^{A'},
\end{split}
\end{gather}
for $Q_G$ and $Q_G^A$, and
\begin{gather}\label{fGA'}
f_G^{A'} =f\left(\frac{[\alpha'g]+0{.}5-pg}{\sqrt{pqg}}\right)
\end{gather}
for $f_G^A$.

We present some examples of applying the above formulas.  For $\mu=-0{.}3$,
$\sigma=10$, $2\beta=0{.}92$, $n=300$ participants, and $s=1000$ steps,
Fig.~2 shows the expected capital increments of the egoists and group
members vs.  the decision threshold, provided that the group makes use of
different voting principles.

\begin{figure}[t]
\centerline{\includegraphics[clip]{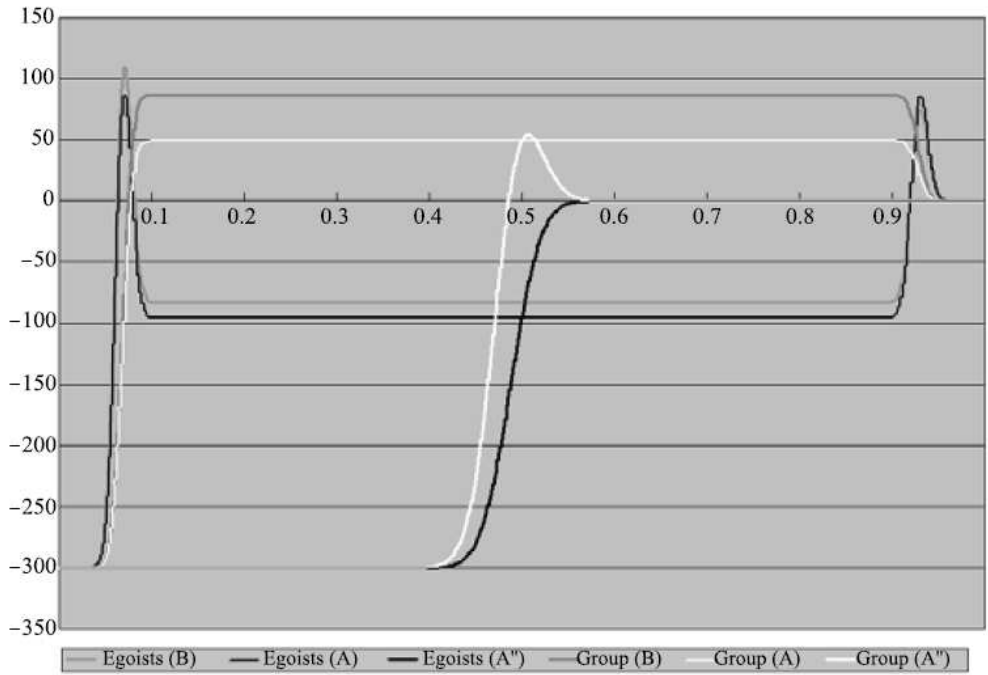}}
\caption{Capital increments of the egoists and group members following
Principles~A, B, and~A$''$ vs.  the decision threshold $\alpha$.  Values of
the parameters:  $\beta=0{.}07$; $n=450$; $\mu=-0{.}3$; $\sigma=10$;
$s=1000$.}
\end{figure}

\begin{figure}[t]
\centerline{\includegraphics[clip]{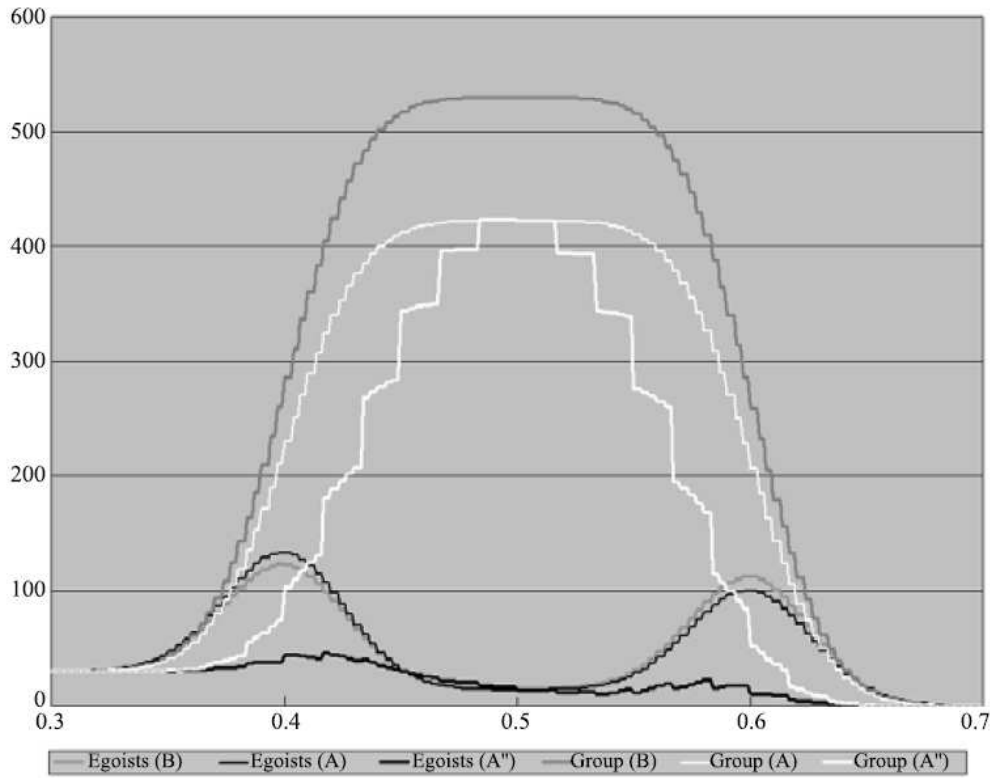}}
\caption{Capital increments of the egoists and group members following
Principles~A, B, and~A$''$ vs.  the decision threshold $\alpha$.  Values of
the parameters:  $\beta=0{.}4$; $n=300$; $\mu=0{.}03$; $\sigma=10$;
$s=1000$.}
\end{figure}

The functions in Fig.~2 are stepwise because the results of voting are
defined by the ratio of the number of participant voting ``for'' and the
{\it integer part} of the threshold $\alpha n$.  The capital increment of a
group member in the case of Principle~A$'$ is shown by the white line which
for $\alpha<0{.}46$ lies slightly below the bold black line corresponding
to the capital increment of a group member in the case of Principle~A.  At
the same time, at the passage to A$'$ the egoists lose more:  for
$\alpha<0{.}46$, the thin black line that represents the variation of their
mean capital passes much lower the bold white line corresponding to
Principle~A.  The capital of a group member grows to the limit for
$\alpha\approx0{.}494$ by 410 units, whereas in the case of accepting all
proposals the reduction in capital is 300 units.  The expected capital
increment of the egoist is positive for $\alpha>0{.}52$ and maximal for
$\alpha\approx0{.}545$.  The group can use a more radical approach to the
choice of the internal threshold assuming that it is equal to the external
one.

{\bf Principle~A$''$}.  {\it The group votes ``for'' a proposal of the
environment if and only if as the result of accepting it the part of its
members getting a positive capital increment exceeds $\alpha$\/}.

In this case, the capital increments are calculated as in the case of
Principle~$'$ with $[\alpha g]$ in place of $[\alpha'g]$.  An example of
using Principle~A$''$ is depicted in Fig.~3, all parameters being the same
as in Fig.~2.

The passage from Principle~A$'$ to Principle~A$''$ in the main worsens the
group activities and, additionally, changes completely the dynamics of the
egoists' capital.  In the example of Fig.~4, practically all proposals of
the unfavorable environment are accepted for $\alpha<0{.}4$, and
practically no proposal is accepted for $\alpha>0{.}57$.

For ${\alpha\!\approx\!0{.}508}$, the expected capital increment of a group
member is maximal:  $M(\widetilde{d}_{\mathcal{G}}^{A''}) \approx54$; the
capital increment of the egoist does not become positive under any
$\alpha$.

The main distinction of Principle~A$''$ lies in that it leads to actual
disappearance of the maxima of capital increments of the egoists, provided
that their fraction is excessive.  In the case of Principles~A and~B, these
peaks are especially high in zones~4 and~5 namely for a small number of
egoists (see Fig.~4), and the passage to Principle~A$''$ eliminates them.
For 300 participants, the height of these peaks is negligible if the
egoists make up less than two thirds of the participants; low peaks for
80\% of egoists (low ``strongly plowed hills'') can be seen in Fig.~5.

If an optimality criterion is defined, then one may raise the question of
the optimal choice of the intragroup voting threshold.  Similarly, it is
possible to vary the threshold of proposal support for principles like B.
In the case of favorable environment, in particular, the requirements on
the mean one-step capital increment of a group member can be set at a
higher level than in an unfavorable environment.  All examples show that to
a certain extent the group can ``dictate its will.'' It deserves noting
that for 300 participants, neutral environment, and simple-majority
procedure, the ratio of the expected capital increment of the member of a
three-participant group to the capital increment of an egoist is
approximately $1{.}4$ in the case of Principle~A and $1{.}75$ in the case
of Principle~B.

\section{conclusions}

Exact and approximate (based on the normal approximation of the binomial
distribution) formulas for the expectations of the capital increment of the
``egoists'' and the ``group members'' were obtained for the considered
model of voting in an environment with certain random parameters.  These
formulas enable one to determine the form of the model trajectory for any
values of the parameters.  The distinctions of different voting principles
were considered using several examples.
\pagebreak

\appendix{}

\PLE{1} Let $(\eta_1,\ldots,\eta_{\ell})=(\zeta_1 I(\overline{\zeta},
\ell_0),\ldots, \zeta_{\ell}I(\overline{\zeta},\ell_0))$.  We denote
$n^+(\overline{\zeta})$ by $n^+$ for brevity.  Using the formula of total
probability for expectation, we get
\begin{gather}\label{Mvot1}
M(\eta_k) =0\cdot P\{n^+ \leq \ell_0\}
+\sum_{x=[\ell_0]+1}^{\ell}M(\eta_k\mid n^+=x)\,P\{n^+=x\}
\end{gather}
for any $k=1,\ldots,\ell$.  Then
\begin{eqnarray}\label{Mvot2}
M(\eta_k\mid n^+=x) &=& M(\eta_k\mid n^+=x,\,\eta_k  >0)\, P\{\eta_k
>0\mid n^+=x\}
\\
 \nonumber         &+& M(\eta_k\mid n^+=x,\,\eta_k\leq0)\, P\{\eta_k\leq0\mid n^+=x\}.
\end{eqnarray}

It follows from independence of the sample elements that
\begin{gather}\label{Meta}
M(\eta_k\mid n^+=x,\,\eta_k>0)=M(\eta_k\mid\eta_k>0)=M(\zeta_k\mid\zeta_k>0).
\end{gather}
This expression is the mean normal random variable, if positive, which is
determined by integration and is as follows:
\begin{gather}\label{Md>0a}
M(\zeta_k\mid\zeta_k>0)
=p^{-1}\frac{1}{\sqrt{2\pi}\sigma}\int\limits_{0}^{\infty}xe^{-\tfrac{(x-\mu)^2}{2\sigma^2}}dx
=\mu+\frac{\sigma f}{p},
\end{gather}
where $f=f(\mu/\sigma)$.  Similarly,
\begin{gather}\label{Md<0a}
M(\eta_k\mid n^+=x,\,\eta_k\leq0)
=M(\zeta_k\mid\zeta_k\leq0) =\mu-\frac{\sigma f}{q}.
\end{gather}
We note that $P\{\eta_k>0\mid n^+=x\}=x/\ell$ and $P\{\eta_k\leq0\mid
n^+=x\}=1-x/\ell$, substitute the determined values in (A.2)
\begin{gather}
\label{Mvot3}
 M(\eta_k\mid n^+=x) =\left(\mu+\frac{\sigma
f}{p}\right)\frac{x}{\ell} +\left(\mu-\frac{\sigma
f}{q}\right)\left(1-\frac{x}{\ell}\right) =\mu-\frac{\sigma
f}{q}+\frac{x\sigma f}{pq\ell}
\end{gather}
and (A.1)
\begin{gather}
\label{Mvot4}
 M(\eta_k)
=\sum_{x=[\ell_0]+1}^{\ell}\left(\mu-\frac{\sigma
f}{q}+\frac{x\sigma f}{pq\ell}\right) \left(\begin{array}{c} \ell
\\ x  \end{array}\right)p^xq^{\ell-x},
\end{gather}
which completes the proof of Lemma.  \hfill$\Box$

\PLE{2} We rearrange (\ref{Mvot}) in
\begin{gather}\label{Mvot'}
M(\eta_k) = \sum_{x=[\ell_0]+1}^{\ell} \left(\begin{array}{c} \ell
\\ x
\end{array}\right)p^xq^{\ell-x} \left(\mu-\frac{\sigma f}{ q}+
          \frac{\sigma f}{pq\ell}\cdot\frac{\sum\limits_{x=[\ell_0]+1}^{\ell}x
          \left(\begin{array}{c} \ell \\ x  \end{array}\right)p^xq^{\ell-x}}
  {\sum\limits_{x=[\ell_0]+1}^{\ell}  \left(\begin{array}{c} \ell \\ x  \end{array}\right)p^xq^{\ell-x}}\right).
\end{gather}
The first sum in the right-hand side is $P\{n^+ > \ell_0\}$; it is
approximated by $F(-\ell'_0)$.  The ratio of the two remaining sums is
$M(n^+\mid n^+ > \ell_0)$.  For $\zeta\sim N(\mu,\sigma^2)$,
\begin{gather}\label{NApp}
M(\zeta\mid\zeta>t)
=\mu+\sigma\frac{f\left(\tfrac{\mu-t}{\sigma}\right)}
             {F\left(\tfrac{\mu-t}{\sigma}\right)}
\end{gather}
(this formula stems from (A.4)).  We use it to approximate $M(n^+\mid n^+ >
\ell_0)$ and obtain
\begin{gather}\label{NApp'}
M(n^+\mid n^+ > \ell_0) \approx
p\ell+\sqrt{pq\ell}\frac{f(\ell'_0)}{F(-\ell'_0)}.
\end{gather}
Substitution of these approximations in (A.8) provides
\begin{eqnarray}\label{Mvot'''}
M(\eta_k) 
&\approx& F(-\ell'_0)\left(\mu-\frac{\sigma f}{ q}+
 \frac{\sigma
 f}{pq\ell}\left(p\ell+\sqrt{pq\ell}\frac{f(\ell'_0)}{F(-\ell'_0)}\right)\right)
 \\
 &=&\mu   F(-\ell'_0) + \frac{\sigma f}{\sqrt{ pq\ell}}f( \ell'_0).
 \nonumber
\end{eqnarray}
\hfill$\Box$

\label{Lastpage}
\label{lastpage}
\revred{F.T.  Aleskerov}
\end{document}